\documentclass[12pt]{article} 
\usepackage{geometry} 
\usepackage{natbib}
\usepackage{color}

\usepackage{indentfirst}
\usepackage{titling}
\usepackage{graphicx} 
\usepackage{float} 
\graphicspath{{/}} 
\usepackage{amsmath}
\usepackage{amssymb}
\usepackage{amsthm}
\usepackage{amsfonts}
\usepackage{bm}
\usepackage{multirow}
\usepackage{booktabs}
\usepackage{appendix}
\usepackage{algorithm}  
\usepackage{algorithmic}  
\pretitle{\begin{center}\Huge\bfseries}
\posttitle{\par\end{center}\vskip 0.5em}
\preauthor{\begin{center}\Large}
\postauthor{\end{center}}
\predate{\par\large\centering}
\postdate{\par}

\title{Global optimization framework for real-time route guidance via variable message sign}
\author{Bai Liu$^a$
\quad Ke Han$^b$\thanks{Corresponding author, email: k.han@imperial.ac.uk}
\quad Jianming Hu$^a$\\\vspace{0.2 in}
$^a$\textit{Department of Automation, Tsinghua University, China}
\\
$^b$\textit{Department of Civil and Environmental Engineering}\\
\textit{Imperial College London, UK}}

\begin{document}
\maketitle
\bibliographystyle{apalike}

\begin{center}
\large{Paper submitted to Transportmatrica A\\
Word count: 5,176 words + 7 figures + 2 tables = 7,426}
\end{center}

\begin{abstract}
Variable message sign (VMS) is an effective traffic management tool for congestion mitigation. The VMS is primarily used as a means of providing factual travel information or genuine route guidance to travelers. However, this may be rendered sub-optimal on a network level by potential network paradoxes and lack of consideration for its cascading effect on the rest of the network. This paper focuses on the design of optimal display strategy of VMS in response to real-time traffic information and its coordination with other intelligent transportation systems such as signal control, in order to explore the full potential of real-time route guidance in combating congestion. We invoke the linear decision rule framework to design the optimal on-line VMS strategy, and test its effectiveness in conjunction with on-line signal control. A simulation case study is conducted on a real-world test network in China, which shows the advantage of the proposed adaptive VMS display strategy over genuine route guidance, as well as its synergies with on-line signal control for congestion mitigation. 
\end{abstract}

\section{Introduction}

Variable message sign (VMS) is an important traffic congestion management tool that provides travel information and route guidance in real time. The impact of VMS on network performance has been assessed in a number of studies \citep{lam1996stochastic, chatterjee2002driver, chen2008effects, shi2009variable}. Although the impact varies in different contexts (urban/extra urban, incident management/congestion mitigation), most of these studies conclude that variable message signs have a great potential to influence route choice, and hence traffic performance.

The majority of studies on VMS tend to focus on the following two aspects. 
\begin{itemize}
\item[(1)] Since the effects of the VMS on network congestion and travel delay are largely dependent on the compliance rate, which is fixed or systematically varied, one stream of research focuses on the estimation of compliance rate for a given scenario using different approaches. Some studies use empirical data to infer drivers' reactions to certain messages displayed, such as those collected from floating cars and loop detectors, surveys and virtual driving simulators \citep{ramsay1997route, chatterjee2002driver, hoye2011evaluation, kattan2011information}. Others tend to systematically derive estimation of the compliance rate based on key factors influencing the compliance rate such as drivers' socio-economic characteristics, gender, familiarity with the network, as well as the clarity and accuracy of the information provided \citep{bonsall1992influence, vaughn1993experimental, wardman1997driver, adler2001investigating, chen2003driver, ben2010road}. 

\item[(2)] Along another line of research, the network effect of the VMS route guidance is investigated in detail.  The impact of VMS on traffic performance, such as congestion, delay, and emission, has been studied in \cite{chatterjee2002driver, chen2008effects, lam1996stochastic, mascia2016impact}. Others have focused on developing models for route choice behavior that account for the learning processes enabled by VMS \citep{dell2009fuzzy, chorus2009traveler, vaughn1993experimental}. 

\end{itemize}

All of the studies mentioned above consider the VMS as a means of providing actual travel information or genuine route guidance to travelers in order to alleviate congestion and reduce travel times. Very few studies have investigated the optimal display strategy in response to real-time traffic information. \cite{peeta2001real} propose a VMS control heuristic framework to comply with the situation when traffic incidents occur. \cite{zuurbier2006generating} designs a generic methodology to generate optimal controlled dynamic prescriptive route guidance, which is disseminated via VMS.

In order to mitigate congestion that occurs on a regular basis, it is crucial to understand the network traffic flow characteristics, which may be affected by a number of factors including origin-destination demand, road condition, signal timing, and driving behavior. While a timely and accurate depiction and dissemination of the system state helps to divert traffic in an efficient way via VMS, potential network paradoxes may suggest sub-optimality of such genuine display strategy. For example, the local improvement of traffic along the corridors affected by the VMS may prove globally deteriorating if the VMS strategy does not take into account its subsequent cascading effect on the rest of the network. In addition, as pointed out by \cite{mascia2016impact}, the effectiveness of VMS can be significantly enhanced via coordination with traffic signal controls. And there is a lack of literature that focuses on designing system-optimal displays strategies, possibly incorporating other intelligent transportation measures such as signal timing, that are full responsive to the traffic condition on the network.

In this paper, we explore the full potential of adaptive, on-line VMS display strategy based on the linear decision rule (LDR) approach typically used for real-time decision making \citep{liu2015data}. Our goal is to design an optimal on-line VMS display strategy that makes system-optimal route guidance in full response to real-time information collected from traffic sensors. A key feature of the proposed strategy is that it reacts to real-time information in a sophisticated way according to the LDR; while the resulting VMS display may not appear intuitive at times, it aims at a system-level optimality in an intelligent way that can be trained off-line. 

\begin{figure}[H]
\setlength{\abovecaptionskip}{3mm}
\setlength{\belowcaptionskip}{0mm}
\center{\includegraphics[width=.6\linewidth]{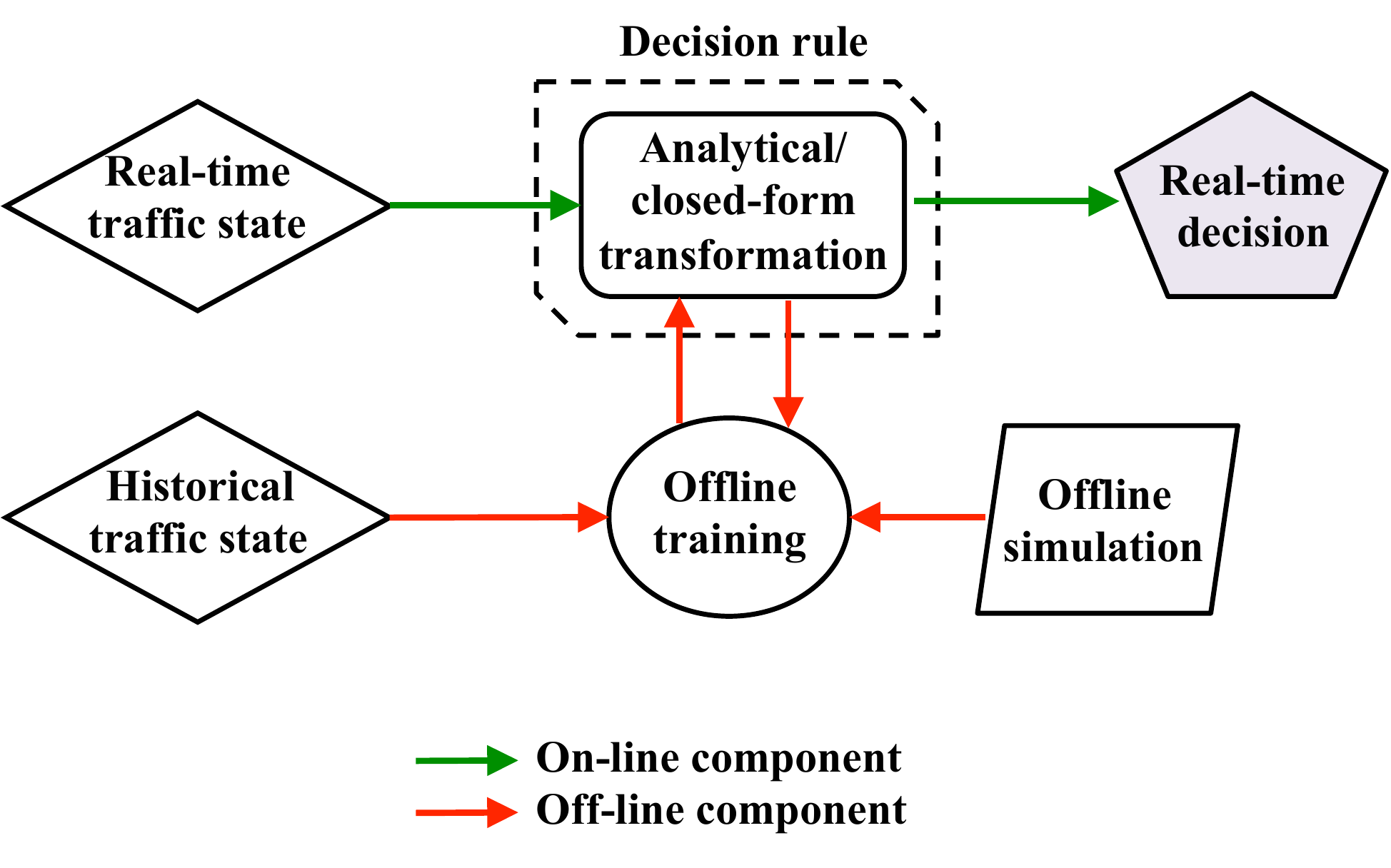}}
\caption{Flow of the decision rule approach.}
\label{figDRflow}
\end{figure}

The LDR framework, illustrated in Figure \ref{figDRflow}, relies on an off-line module and an on-line module. The on-line module converts real-time system state to control parameters via analytic or closed-form transformation (or, in the case of LDR, linear transformation). This allows a timely decision making in a real-time environment. The performance of the such a transformation, on the other hand, needs to be constantly improved using the off-line training module, which takes into account historical traffic information. 

The LDR-based VMS display strategy reactions timely to the network condition, in ways that guarantee system-wide optimal due to the extensive off-line training. Moreover, this approach allows the objective of the on-line control to be arbitrarily defined, as long as it can be computed via the off-line simulation. Last but not least, as \cite{liu2015data} have demonstrated, the LDR approach can be readily applied to responsive signal control, making possible the coordination between VMS and signal control, which will be explored in this paper.



The proposed methodology is applied to a real-world test site in the city of Haining, China, where a VMS is installed at a key junction for diverting traffic during afternoon peak hours.  Traffic counts collected at different locations for 20 days in the test network are used to perform the off-line training, i.e. to find the optimal linear coefficients in the LDR. The proposed adaptive VMS strategy is compared with the more naive strategy that shows the accurate travel information, and is demonstrated to yield improved system performance. Sensitivity analysis is also conducted to verify the effectiveness of the adaptive strategy against different levels of compliance rates.

The rest of this paper is organized as follows. Section \ref{secDR} presents the framework of on-line decision making based on the decision rule approach. We then discuss a real-world case study of VMS control and the implementation of a linear decision rule in Section \ref{secCase}. We show the numerical test results for this case study in Section \ref{secNumerical}. Finally, Section \ref{secConclusion} offers some concluding remarks.

\section{Decision Rule for Real-Time Decision Making}\label{secDR}

\subsection{The conceptual model}

We begin with illustrating the general framework of decision rule and its application to real-world traffic management, without referring to specific network structure, data type, or control mechanism. Mathematically, let us denote by $\boldsymbol q$ a vector of state variables that describe the traffic system in the present and/or a few moments in the past. Generally speaking $\boldsymbol q$ may represent arbitrary traffic quantities such as flow, density, speed, or any analytical transformation of them through data fusion or system reconstruction techniques. The decision rule is instantiated as the following analytical mapping:
\begin{equation}\label{eqnDR}
u~=~P_{\Omega}\big[ f(x,\,\boldsymbol q)\big]
\end{equation}
where $x$ consists of parameters to be determined via the offline training, and $u$ is the set of feasible controls to be implemented. The functional form of $f(\cdot)$ is arbitrary, and can be linear, nonlinear, or implicitly defined (e.g. artificial neural networks). $P_{\Omega}[\cdot]$ denotes the projection onto the set of feasible control parameters $\Omega$, which may be characterized by linear or nonlinear constraints in continuous or discrete domains. For example, in the case of VMS the set $\Omega$ may contain finite number of routing messages, and is thus discrete. A decision rule like \eqref{eqnDR} can yield timely decision that allows real-time operation as it involves analytical/closed-form transformation, and its efficacy can be improved to a target degree of optimality through the off-line training of the parameters $x$.

Let $\Phi(\boldsymbol q,\,u)=\Phi\left(\boldsymbol q,\, P_{\Omega}\big[f(x,\,\boldsymbol q)\big]\right)$ be an arbitrary network performance function, which obviously depends on the system state $\boldsymbol q$ and the control $u$ (along with some other uncertain or endogenous variables). For example, $\Phi$ may be the delay at a particular junction, or the total emission along certain corridor. Without loss of generality we assume that $\Phi$ is subject to minimization in this paper. Therefore, the problem of optimal decision rule can be formulated as:
\begin{equation}\label{eqn2}
\min_x \Phi\left(\boldsymbol q,\, P_{\Omega}\big[f(x,\,\boldsymbol q)\big]\right)
\end{equation}
\noindent Note that, however, $\boldsymbol q$ is a stochastic variable that changes on a daily basis. For example, $\boldsymbol q$ can be the vector of time-varying on-ramp demands in a highway network, which vary from day to day. Therefore, a robust feedback control policy such as \eqref{eqnDR} must take into account the stochasticity in the system state $\boldsymbol q$. With this in mind, the off-line training of the decision rule may be conceptually formulated as the following stochastic optimization problem
\begin{equation}\label{eqn3}
\min_x \mathbb{E}_{\mathbb{D}}\left(\boldsymbol q,\, P_{\Omega}\big[f(x,\,\boldsymbol q)\big]\right)
\end{equation}
\noindent where $\mathbb{D}$ characterizes the stochasticity of $\boldsymbol q$; e.g. it can be the joint distribution of the individual components of $\boldsymbol q$. Note that in order to solve \eqref{eqn3}, we need to access a set of historical data on the uncertain variable $\boldsymbol q$, as shown in Figure \ref{figDRflow}.

\subsection{Solution method}
The stochastic optimization problem \eqref{eqn3} is difficult to solve for the following reasons:
\begin{itemize}
\item[(1)] The objective function $\Phi$ is likely to be highly nonlinear and non-convex, making the stochastic optimization problem intractable. 
\item[(2)] The stochastic vector $\boldsymbol q$ is usually high-dimensional, and the joint distribution of its components is very difficult to estimate in practice. 
\end{itemize}

We propose a practical way of solving the problem \eqref{eqn3} using a Monte-Carlo approach. To fix the idea, we treat the performance function $\Phi\left(\boldsymbol q,\, P_{\Omega}\big[f(x,\,\boldsymbol q)\big]\right)$ as a single-valued random variable parameterized by $x$; its randomness is due to the random state variable $\boldsymbol q$ (and possibly some other random variables endogenous to the system). The corresponding stochastic optimization problem becomes
\begin{equation}\label{DROMC}
\min_x  \mathbb{E}_{\mathbb{D}(x)}\Phi\left(\boldsymbol q,\, P_{\Omega}\big[f(x,\,\boldsymbol q)\big]\right)
\end{equation}
\noindent $\mathbb{D}(x)$ is the distribution of the $x$-dependent random variable $\Phi\left(\boldsymbol q,\, P_{\Omega}\big[f(x,\,\boldsymbol q)\big]\right)$.

Consider a set of $K$ sampled historical data of the system state variable: $\{ \boldsymbol q^{(1)},\ldots, \boldsymbol q^{(K)}\}$. For each given $x$, we can obtain a sequence of objective values 
$$
\boldsymbol \Phi^i(x)\doteq \Phi\left(\boldsymbol q^{(i)},\, P_{\Omega}\big[f(x,\, \boldsymbol q^{(i)})\big]\right),\qquad 1\leq i\leq K,
$$
which are $K$ independent samples  drawn from $\mathbb{D}(x)$. Assuming sufficient sample size $K$, we may then approximate the expectation in \eqref{DROMC} using the average of $\boldsymbol\Phi^i(x)$, $1\leq i\leq K$.

This Monte-Carlo type procedure calculates the objective function of the optimization problem \eqref{DROMC} for a given $x$. Therefore, this procedure can be easily integrated with metaheuristic optimization methods that can search for the optimal $x$ while maintaining an efficient trade-off between system optimality and computational overheads. For our particular application presented below, we employ the {\it particle swarm optimization} \citep{banks2007review}, while noting that other metaheuristic methods such as simulated annealing and genetic algorithm can be equally applied.

\section{Case Study}\label{secCase}

\begin{figure}[H]
\setlength{\abovecaptionskip}{3mm}
\setlength{\belowcaptionskip}{0mm}
\center{\includegraphics[width=\linewidth]{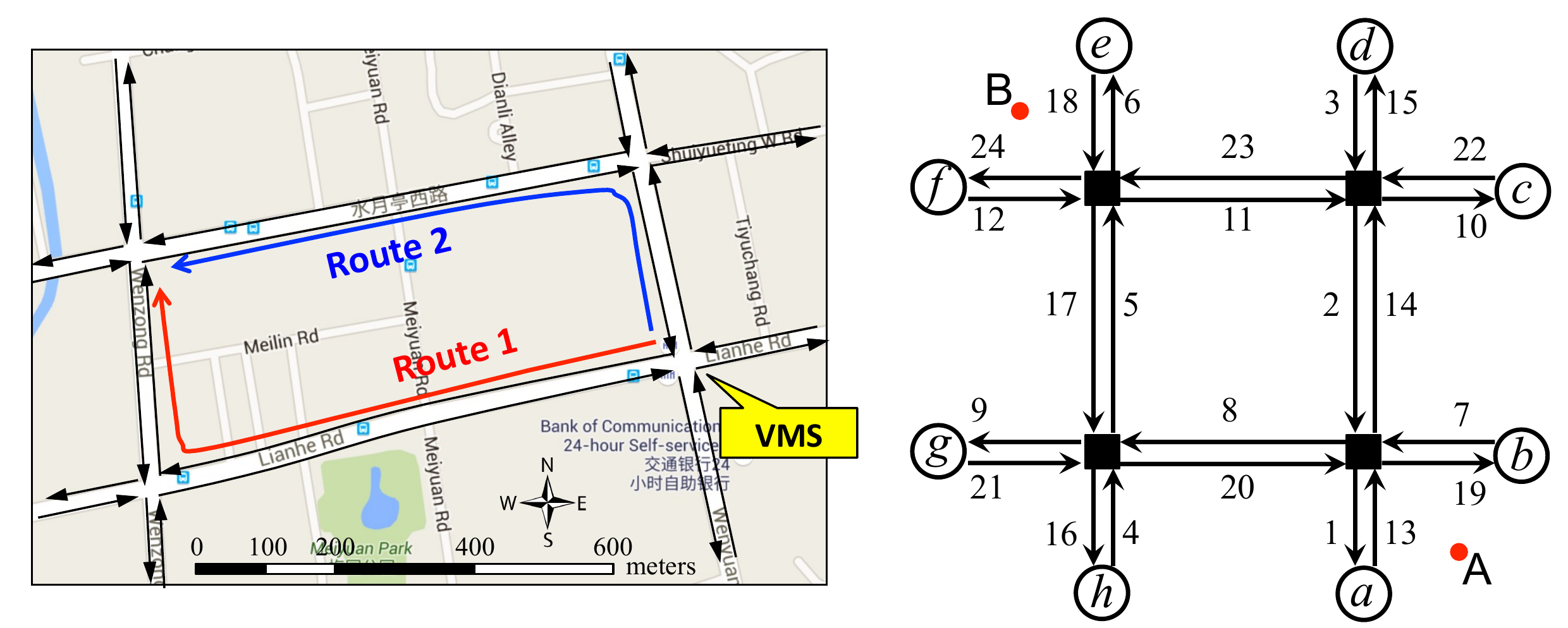}}
\caption{The test network in Haining, China}
\label{figVMS_China}
\end{figure}

A case study of the proposed real-time adaptive VMS strategy is carried out in the city of Haining, China. Figure \ref{figVMS_China} shows the layout of the test network, which consists of 24 links and 4 intersections. 

We focus on the afternoon peak from 4 - 5pm, when the network is dominated by traffic traveling towards Northwest. For drivers traveling from point $A$ to point $B$ (see Figure \ref{figVMS_China}), two alternative routes are available as shown in the figure. The two routes have similar travel distance and level of congestion, and are both viable options for the drivers. The location of the VMS is highlighted, which provides critical route guidance for drivers approaching the network from the Southeast (location $A$).

\subsection{Modeling scenarios and assumptions}\label{subsecMSAA}

The VMS is located at point $A$ (see Figure \ref{figVMS_China}). The following scenarios/assumptions are considered regarding the effect of the route guidance provided by the VMS.
\begin{itemize}
\item The VMS is only visible to traffic approaching from links $7$ and $13$. Moreover, the message will only influence drivers belonging to origin-destination (OD) pairs $(a,\,e)$, $(a,\,f)$, $(b,\,e)$, $(b,\,f)$. 

\item Five types of messages can be displayed by the VMS: (1) Route 1 -- strongly; (2) Route 1 -- moderately; (3) no display (or neutral suggestion); (4) Route 2 -- moderately; (5) Route 2 -- strongly. Each display corresponds to a certain compliance rate. 

\item The rest of the traffic not influenced by the VMS are called {\it background traffic}. The demand profile and route choices of background traffic are extracted from empirical data collected on the test site. 
\end{itemize}

We associate the VMS display with the compliance rate in the following way. We let $c_i\in (0,\,1)$ be the proportion of drivers from zones $a$ or $b$ who choose Route 1, when the VMS displays strategy ($i$), $i=1,\ldots, 5$. And we use $\mathbf{P}_{route 1}=(c_1,\,\ldots,\,c_5)$ to denote a complete set of compliance rates. Later in our numerical study, we will vary the choices of $c_i$ in order to test the sensitivity of the proposed control to the compliance rates, which are in general very difficult to estimate in a real-world traffic network.

Besides the VMS, we also consider traffic signal coordination at the relevant junctions in order to maximize the effect of route guidance provided by the VMS. The design and optimization of the on-line signal timing based on the linear decision rule are detailed in \cite{liu2015data} and omitted here.

\subsection{Agent-based simulation}\label{subsecABS}

We employ an agent-based model in order to capture the individual routing decision influenced by the VMS. Each agent has his/her own origin and destination, and is characterized as (a) background agent whose routing decision is not affected by the VMS; and (b) target agent, whose routing is influenced by the VMS. In this case study, the target agents are those associated with the OD pairs $(a,\,e)$, $(a,\,f)$, $(b,\,e)$ and $(b,\,f)$. For the propagation of link flow, we employ a link performance function \citep{friesz1993variational}, which specifies the link travel time as a function of link occupancy:
\begin{equation}\label{LDM}
D(t)~=~F(X(t))~=~\alpha X(t) +\beta
\end{equation}
where $D(t)$ denotes the travel time on the link when the time of entry is $t$; $X(t)$ is the link occupancy at the entry time $t$; $\alpha$ and $\beta$ are positive constants. \\

\noindent{\bf Remark.}
{\it In this paper we use the simple link performance function \eqref{LDM} for simplicity. The proposed methodological framework, however, allows easy extension of the network simulation model to incorporate more realistic features of a dynamic network such as car following and  spillback. In particular, analytical traffic network models (e.g. kinematic wave models) and microsimulation models can be integrated with little modification of the method.}

Regarding vehicle movement at the intersections, we distinguish between background agents and target agents. For the background agents, we assume that their movements at intersections are prescribed by turning probabilities derived from historical data collected on the test site. For the target agents, their routing choices are explicitly enumerated and hence determine the intersection movements. Note that the routing of target agents are also affected by the VMS with certain probability (compliance rate), which is defined at the end of Section \ref{subsecMSAA}.

Due to the randomness in the routing of both background and target agents, the proposed agent-based simulation is inherently stochastic. Note that such systematic stochasticity, together with the random state variable $\boldsymbol q$, are manifested in the objective variable $\Phi\left(\boldsymbol q,\, P_{\Omega}\big[f(x,\,\boldsymbol q)\big]\right)$; see \eqref{eqn3}. The proposed Monte-Carlo approach is able to encapsulate these various sources of randomness, and support the optimization \eqref{eqn3} via metaheuristic solution methods.

\subsection{Implementation of the linear decision rule approach}

We consider the linear decision rule, which specifies the response (control $u$) as a linear transformation of the observation (system state $\boldsymbol q$). Mathematically, we write the state variable as
$$
\boldsymbol q~=~\left(q(1),\, q(2),\,\ldots,\, q(\mathcal{T})\right)
$$
where each $q(t)$ represents a vector of system state at time step $t$, $t=1,\ldots, \mathcal{T}$. In our case study, 
$$
q(t)~=~\left(X_1(t),\,  X_2(t),\,\ldots,\, X_{24}(t)\right)^T
$$ 
where $X_i(t)$ is the average link occupancy on link $i$ during the time step $t$. At each time step $t$, the linear decision rule specifies that
\begin{equation}\label{LDRcase}
\mu(t)~=~f(\mathbf{A},\,\boldsymbol q)~=~\mathbf{A} \cdot 
\begin{bmatrix}
q(t-1)
\\
q(t-2)
\\
\ldots
\\
q(t-\delta)
\end{bmatrix}\qquad 1\leq \delta \leq t-1
\end{equation}
where $\mathbf{A}$ (equivalent to $x$ in Eqn \eqref{eqnDR}) is a row vector of proper dimension, whose elements are to be optimized via the off-line training phase, i.e. \eqref{eqn3}). Eqn \eqref{LDRcase} implies that the control decision made for time step $t$ may depend on the system state in the previous $\delta$ time steps. According to \eqref{eqnDR}, we need to project $\mu(t)$ onto the set $\Omega$ of feasible controls, which in this case are the five VMS messages:
\begin{equation}
u(t)~=~P_{\Omega}[\mu(t)]~=~
\begin{cases}
\hbox{Route 1 -- strongly} \qquad &\mu(t) \in (-\infty,\, m_1]
\\
\hbox{Route 1 -- moderately} \qquad & \mu(t) \in (m_1,\, m_2]
\\
\hbox{No display} \qquad & \mu(t)\in (m_2,\,m_3]
\\
\hbox{Route 2 -- moderately} \qquad & \mu(t) \in (m_3,\, m_4]
\\
\hbox{Route 2 -- strongly} \qquad & \mu(t) \in (m_4,\, \infty)
\end{cases}
\end{equation}
Here, we are using a simple linear projection onto five display strategies. The choices of the values $m_1$ - $m_4$ are arbitrary, but may influence the performance of the linear decision rule or the off-line training procedure.

\subsection{Off-line training}\label{subsecofflinetraining}

In order to capture the variability in the daily demand profile, we collect historical link occupancy data during 4 - 5pm from 1 - 10 February 2016 as the training input. The link occupancy data were recorded by microwave detectors and provided by the transportation authority in Haining city, Jiangsu Province. These 10 days of link occupancy data will be used to derive OD demand profiles and estimate vehicle turning probabilities at intersections. These data will also be used in the network loading combined with agent-based simulation discussed in Section \ref{subsecABS}.

The objective function $\Phi$ (see Eqn \eqref{eqn3}) is defined to be the network-wide average travel time per vehicle, which is obtained from the agent-based simulation. We note that the proposed framework and metaheuristic solution method allow much flexibility in the objective function, ranging from speed, throughput, delay, to emission and fuel consumption, measured at different spatial references (e.g. junction level, corridor level, or network-wide level).

For the metaheuristic search method, we employ particle swarm optimization (PSO) \citep{banks2007review}, which offers an efficient and flexible trade-off between optimality of the solution and computational overhead. PSO is based on the
social behavior of a group of animals, called a {\it swarm}. In a swarm, the animals are
represented as particles, and can collaborate and share information to
adjust their positions in the search for a certain location. The
adjustment of their positions is based on the swarm's collective
memory on the best location attained so far ($\mathbf{G}$), and the individual memory of the best location that the
individual particle has attained so far ($\mathbf{P}_{j}$ for particle $j$). As
a result of the position adjustment, the particles tend to converge to
either $\mathbf{G}$ or $\mathbf{P}_j$. Although the performance of PSO varies depending on the domain of application or parameters chosen, research shows evidences
of PSO or its variations outperforming well-established metaheuristics
such as genetic algorithm, ant colony optimization, simulated
annealing, and tabu search.

Given the objective function $f(\cdot)$ to be minimized, Algorithm \ref{alg4.1} summarizes the PSO procedure.

\begin{algorithm}  
\caption{Particle swarm optimization}  
\label{alg4.1}  
\begin{algorithmic}   

\FOR {i = 1, 2, ..., iteration}
	\STATE{
	Initialize $\mathbf{X_{i}}$\\
	Initialize $\mathbf{v_{i}}$
	}
\ENDFOR
\FOR {i = 1, 2, ..., iteration}
	\FOR {j = 1, 2, ..., particle}
	\STATE{
		Update the particle's position: $\mathbf{X}_{j}$ $\gets$ $\mathbf{X}_{j}$ + $\mathbf{v}_{j}$\;
		\IF {f($\mathbf{X}_{j}$) $<$ f($\mathbf{P}_{j}$)}
			\STATE Update the particle's best known position: $\mathbf{P}_{j}$ $\gets$ $\mathbf{X}_{j}$\;
			\IF {f($\mathbf{P}_{j}$) $<$ f($\mathbf{G}$)}
				\STATE Update the swarm's best known position: $\mathbf{G}$ $\gets$ $\mathbf{P}_{j}$\;
			\ENDIF
		\ENDIF \\
		Pick random numbers: $r_{1}, r_{2} \sim U(0,1)$\; \\
		Update the particle's velocity: $\mathbf{v}_{j}$ $\gets$ $\mathbf{v}_{j}$ + $c_{1} r_{1}$($\mathbf{P}_{j}$ - $\mathbf{X}_{j}$) + $c_{2} r_{2}$($\mathbf{G}$ - $\mathbf{X}_{j}$)\;}
	\ENDFOR
\ENDFOR
\end{algorithmic}  
\end{algorithm} 

\noindent In our case study, we set the number of particles to be 30 and the number of major iterations in the PSO algorithm to be 20.

\section{Numerical Study}\label{secNumerical}

\subsection{Description of the data}

We obtain real data from the transportation authority in Haining city, Jiangsu Province in China. The flow data is collected by microwave detectors and covers the12 roads in both directions (see Figure \ref{figVMS_China}). The time resolution of the data is 1 minute, and the recording period spans are 16:00-17:00 from 1 to 20 February, 2016. 

\begin{figure}[H]
\setlength{\abovecaptionskip}{3mm}
\setlength{\belowcaptionskip}{0mm}
\center{\includegraphics[width=.9\textwidth]{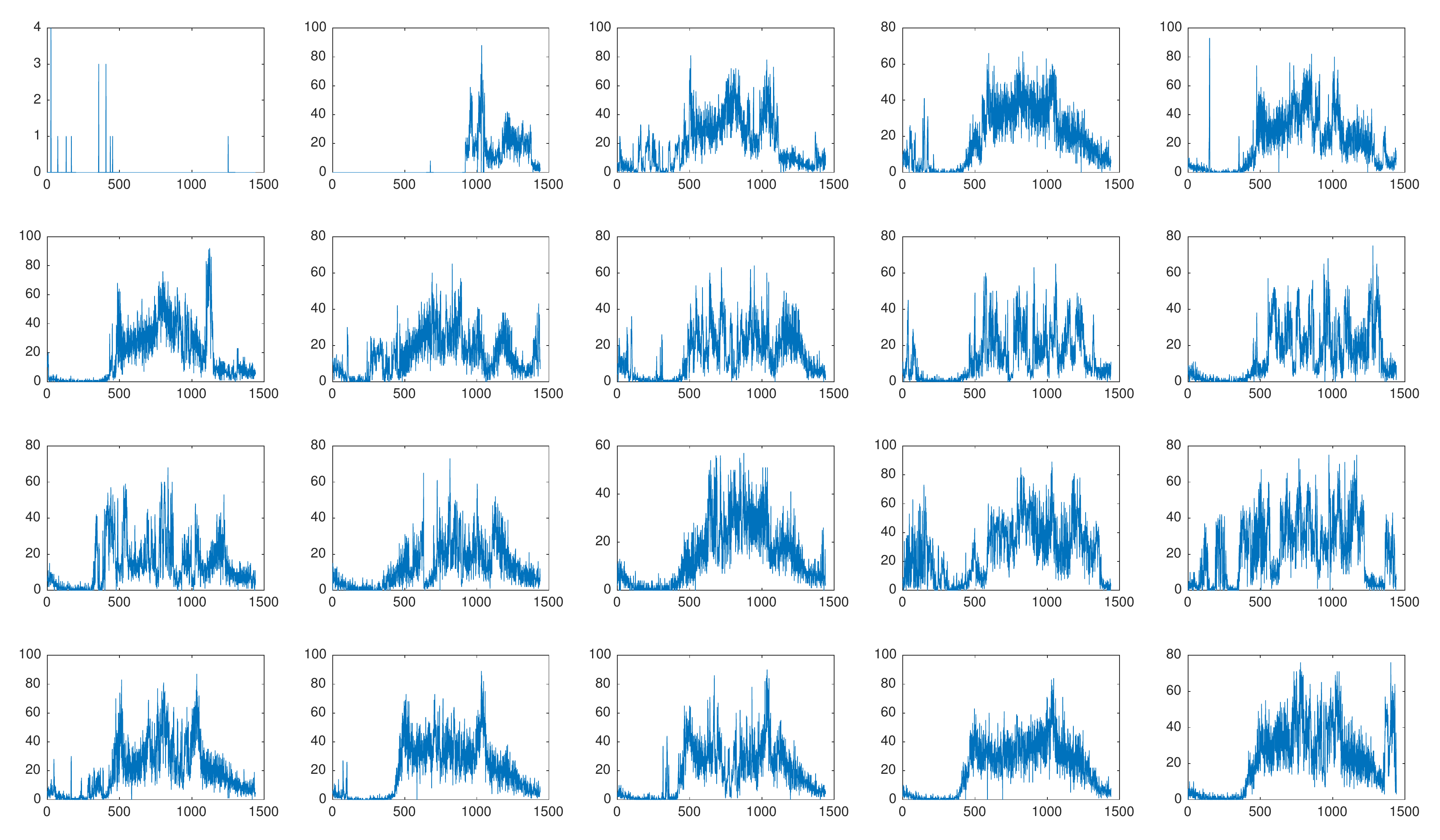}}
\caption{Flow data for 20 different days on one road. x-axis: time (in min); y-axis: flow (in vehicle/min)}
\label{figdatasample}
\end{figure}

One sample dataset is shown in Figure \ref{figdatasample} corresponding to the 20 days of flow data on the road connecting zone $A$ with the rest of the network (see Figure \ref{figVMS_China}). 

As can be seen from Figure \ref{figdatasample} that the raw dataset has some missing/incomplete data. In addition, the road-level flow data is not distinguished by the travel direction. Thus additional data extrapolation techniques are used to reconstruct the dataset for the use of the decision rule approach, for which the details are not elaborated due to space limitation.

%

\subsection{Simulation-based optimization}
As we present in Section \ref{subsecofflinetraining}, the off-line training of the linear decision rule boils down to a stochastic optimization, which requires iteratively perform the agent-based simulation with given training data. We use the data from 1 to 10 February as the training data for the off-line computation, and the rest 10 days of data for testing. 

We use the particle swarm optimization to perform the stochastic optimization, with the following optimization parameters:

\begin{table}[H]
\setlength{\abovecaptionskip}{0mm}
\setlength{\belowcaptionskip}{3mm}
\centering
\caption{PSO parameters}
\label{Tab:OptPara}
 \begin{tabular}{cc}
  \toprule
  \textbf{Parameter} & \textbf{Value} \\
  \midrule
\# of particles & 20 \\
 max \# of iterations & 30 \\
  \bottomrule
\end{tabular}
\end{table}

We use the data in 16:00-17:00 from Feb.11 to Feb.20 as testing data for the linear decision rule approach. For comparison purposes, the following four alternatives are considered to extensively verify the effectiveness of our approach. 

\begin{itemize}
\item \textbf{[Genuine VMS + default signal].} The VMS display reflects the actual traffic conditions on the two alternatives, and the traffic signals have the default timing plan.
\item \textbf{[Genuine VMS + coordinated signal].} The VMS display reflects the actual traffic conditions on the two alternatives, and the traffic signal control follow an LDR-based responsive control
\item \textbf{[LDR VMS + default signal].} The VMS display follows the LDR-based responsive control, while the traffic signals have the default timing plan.
\item \textbf{[LDR VMS + coordinated signal].} Both VMS display and traffic signal control follow LDR-based responsive control
\end{itemize}
\noindent Here, the default signal timing is set by assigning equal green time to the four signal phases at the intersections. The coordinated signal is obtained through the LDR approach, and the optimal LDR for the signal control is trained together with the LDR for the VMS.

We simulate under the situation that the number of objective agents is 80\% of background vehicles. We also choose 3 different sets of compliance rates $\mathbf{P}_{route 1}$ to examine the sensitivity and robustness of our proposed VMS control strategy. These are compliance rates are: $\mathbf{P}_{route 1}=(0.3, 0.4, 0.5, 0.6, 0.7)$, $(0.2, 0.4, 0.5, 0.6, 0.8)$ or $(0.1, 0.3, 0.5, 0.7, 0.9)$. We refer the reader to the end of Section \ref{subsecMSAA} for an explanation of these compliance rates.

\subsection{Results and analysis}
\begin{figure}[H]
\setlength{\abovecaptionskip}{3mm}
\setlength{\belowcaptionskip}{0mm}
\center{\includegraphics[width=0.75\linewidth]{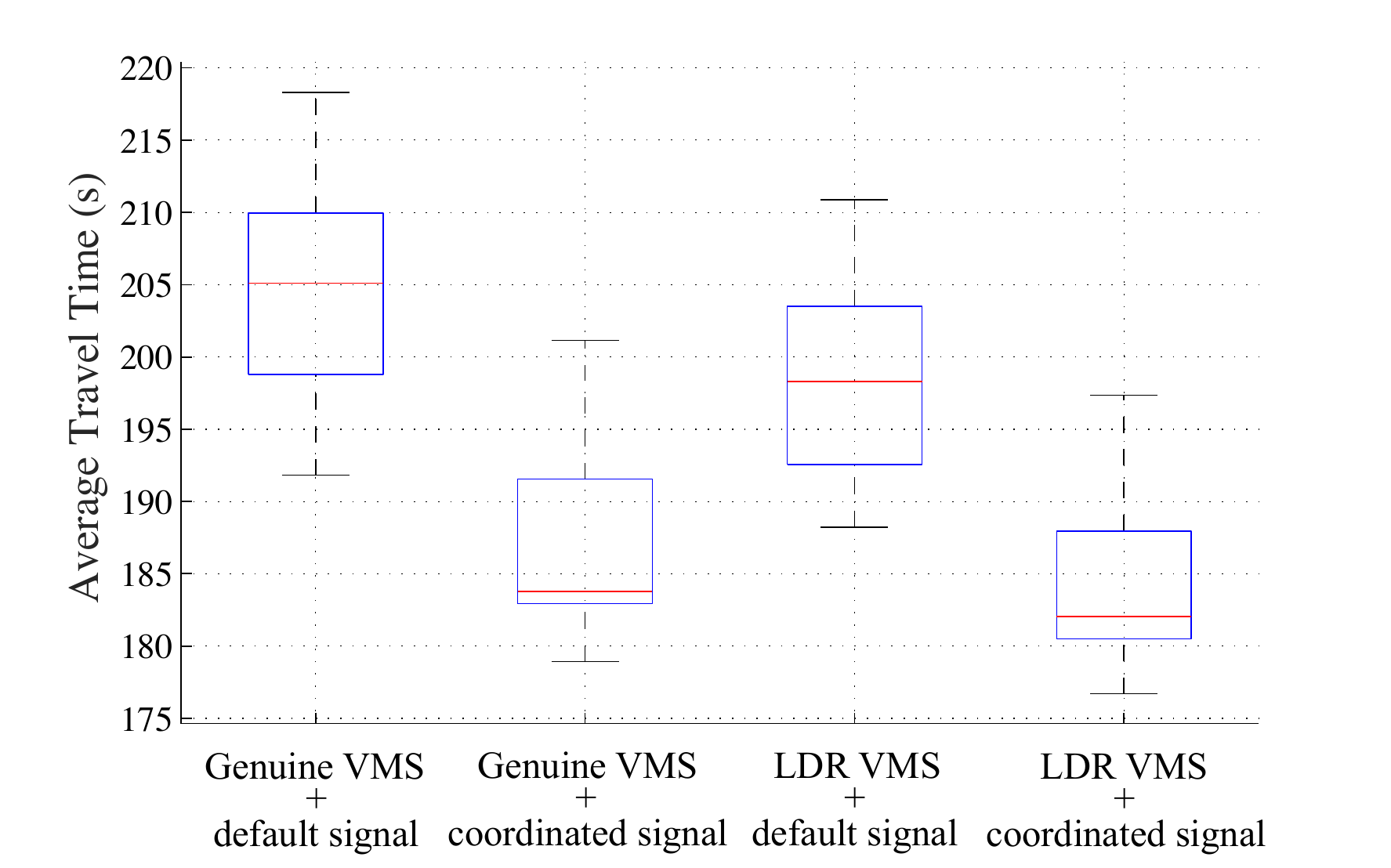}}
\caption{Comparison of Different Guidance Policies ($\mathbf{P}_{route 1} = (0.3, 0.4, 0.5, 0.6, 0.7)$)}
\label{Fig:Result0607}
\end{figure}

\begin{figure}[H]
\setlength{\abovecaptionskip}{3mm}
\setlength{\belowcaptionskip}{0mm}
\center{\includegraphics[width=0.75\linewidth]{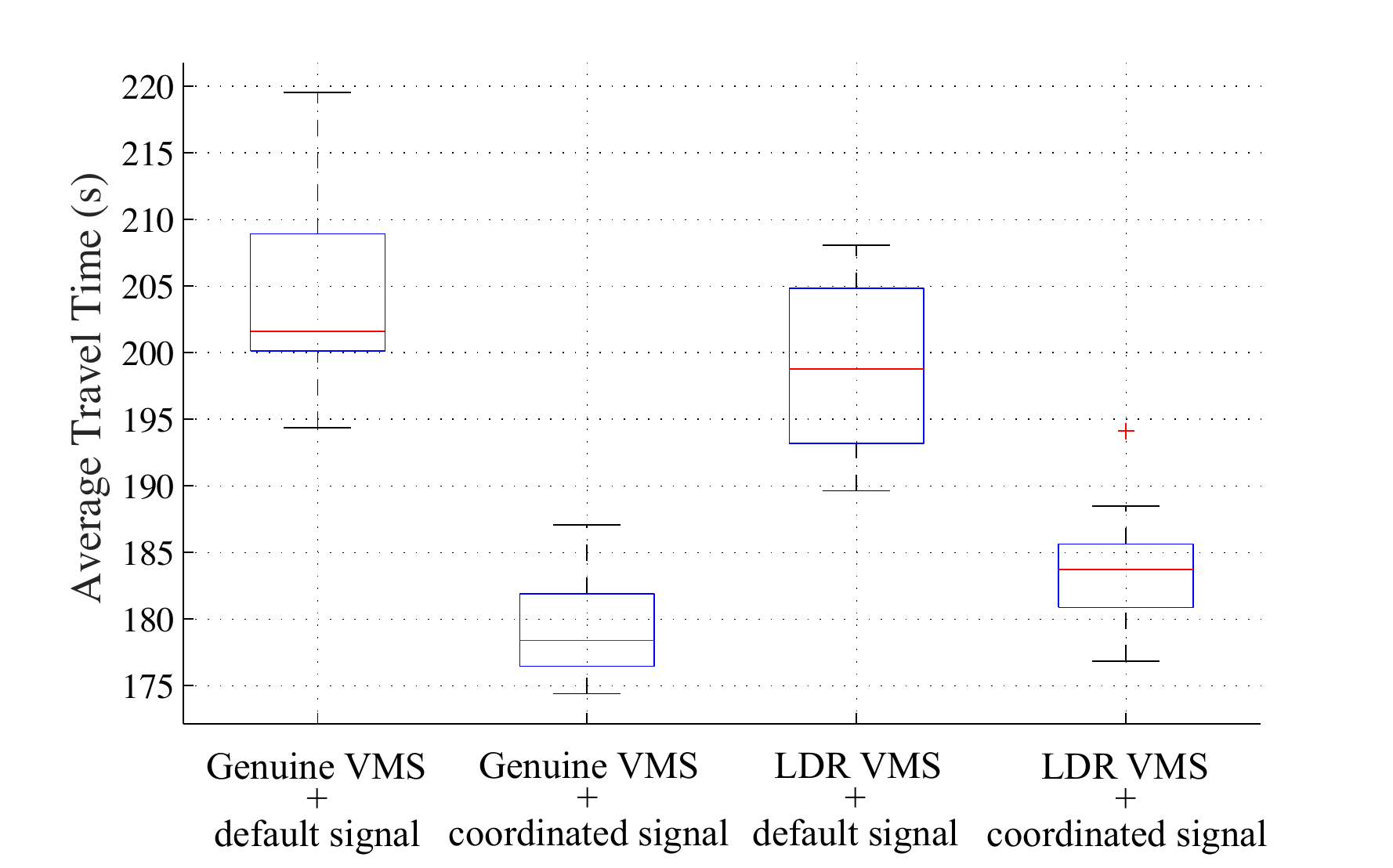}}
\caption{Comparison of Different Guidance Policies ($\mathbf{P}_{route 1} = (0.2, 0.4, 0.5, 0.6, 0.8)$)}
\label{Fig:Result0608}
\end{figure}

\begin{figure}[H]
\setlength{\abovecaptionskip}{3mm}
\setlength{\belowcaptionskip}{0mm}
\center{\includegraphics[width=0.75\linewidth]{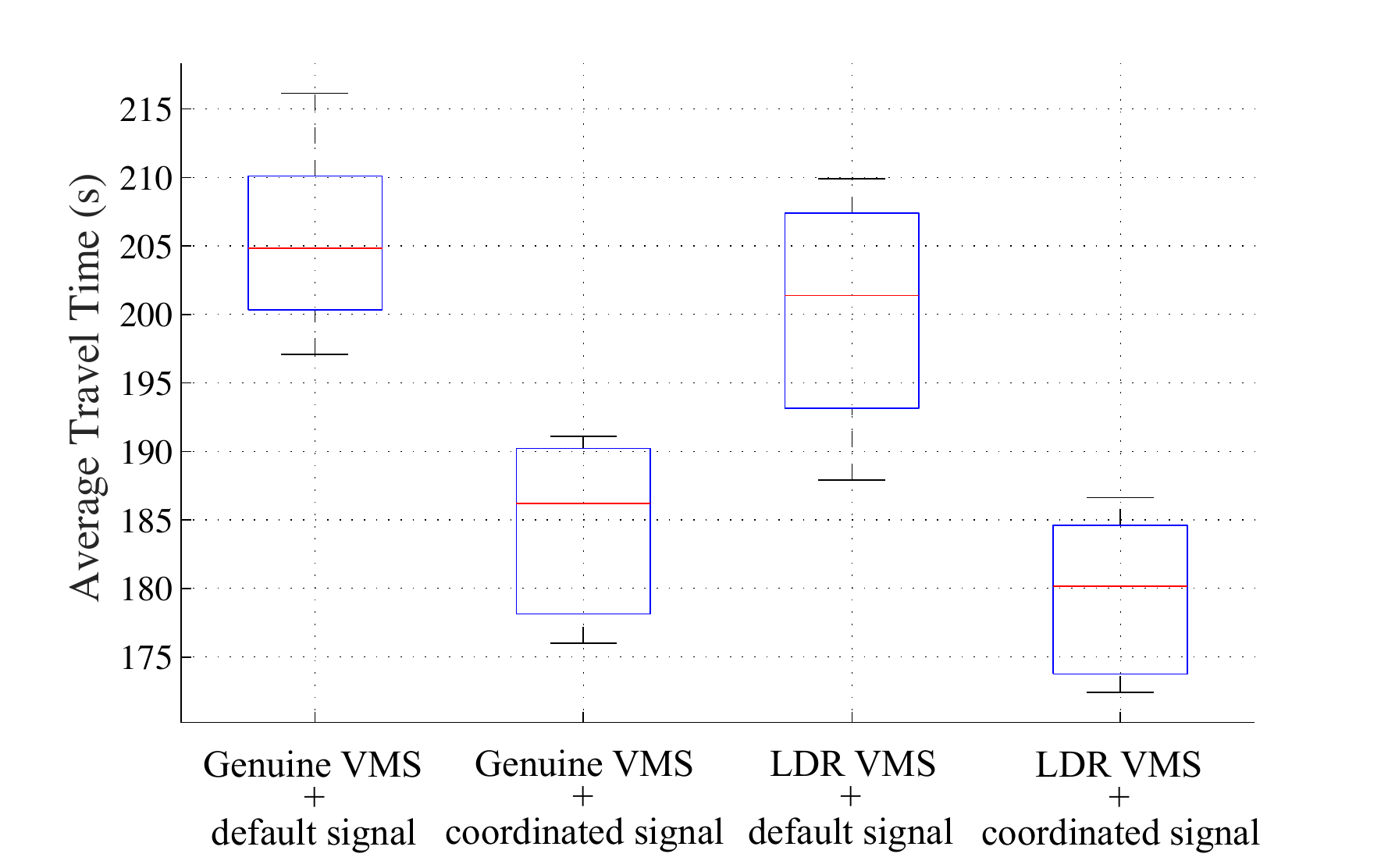}}
\caption{Comparison of Different Guidance Policies ($\mathbf{P}_{route 1} = (0.1, 0.3, 0.5, 0.7, 0.9)$)}
\label{Fig:Result0709}
\end{figure}

For each choice of $\mathbf{P}_{route 1}$, we test the four different control strategies based on the test dataset. The performances of the four strategies in terms of the average travel time per vehicle in the network are shown in Figures \ref{Fig:Result0607}, \ref{Fig:Result0608}, and \ref{Fig:Result0709} for $\mathbf{P}_{route 1} = (0.3, 0.4, 0.5, 0.6, 0.7)$, $(0.2, 0.4, 0.5, 0.6, 0.8)$, and $(0.1, 0.3, 0.5, 0.7, 0.9)$, respectively.

\noindent In addition, the mean objectives (average travel time per vehicle) over the simulation test runs are summarized in Table \ref{Tab:AverageTime}.

\begin{table}[H]
\setlength{\abovecaptionskip}{0mm}
\setlength{\belowcaptionskip}{3mm}
\centering
\caption{Mean objectives (average travel time per vehicle) under different choices of the compliance rates $\mathbf{P}_{route 1}$}
\label{Tab:AverageTime}
 \begin{tabular}{ccccc}
  \toprule
  $\mathbf{P}_{route 1}$ & $\mathbf{(0.3, 0.4, 0.5, 0.6, 0.7)}$ & $\mathbf{(0.2, 0.4, 0.5, 0.6, 0.8)}$ & $\mathbf{(0.1, 0.3, 0.5, 0.7, 0.9)}$\\
  \midrule
  Genuine VMS +  & $204.38$ s & $204.95$ s & $205.75$ s \\
  default signal  & & & 
  \\\hline
  Genuine VMS +  & $186.42$ s & $179.35$ s & $184.69$ s \\
  coordinated signal & & & 
  \\\hline
  Optimal VMS +   & $198.11$ s & $198.56$ s & $200.07$ s \\
  default signal  & & & 
  \\\hline
  LDR VMS +   & $183.96$ s & $184.05$ s & $179.65$ s \\
  coordinated signal & & & \\
  \bottomrule
\end{tabular}
\end{table}

The results show that the LDR-based VMS strategy outperform the genuine display in most of the test cases, with only one exception: Genuine VMS + coordinated signal vs. LDR VMS + coordinated signal, $\mathbf{P}_{route1}=(0.2, 0.4, 0.5, 0.6, 0.8)$. This exception may be due to the probabilistic search performed by the particle swarm optimization, which does not guarantee global optimality of the optimization problem. On average, the trained on-line VMS display saves 3-5 seconds of travel time per vehicle without or without coordination with signals. 

We also observe that incorporating on-line signal control significantly improves the performance of the VMS. However, such substantial margin of improvement is caused by the following factors.
\begin{itemize}
\item[(1)] The naive default signal setting is far from being optimal. On-going work is underway to obtain more realistic default signal timings either from the transportation authority or via simulation-based optimization. 

\item[(2)] The traffic signal controls naturally have a greater influence than the VMS on the performance of the entire test network, since they have total control over all the four intersections, and impact both the subject agents and the background agents. 

\end{itemize}

\noindent We also observe from Table \ref{Tab:AverageTime} that the effectiveness of the VMS is not very sensitive to the different compliance rates. That is, across all three choices of $\mathbf{P}_{route1}$, the saving in travel time of the proposed VMS strategy remains approximately 5-6 s with default signal timing, and 3-5 s with coordinated signal control. 

\begin{figure}[H]
\setlength{\abovecaptionskip}{3mm}
\setlength{\belowcaptionskip}{0mm}
\center{\includegraphics[width=\linewidth]{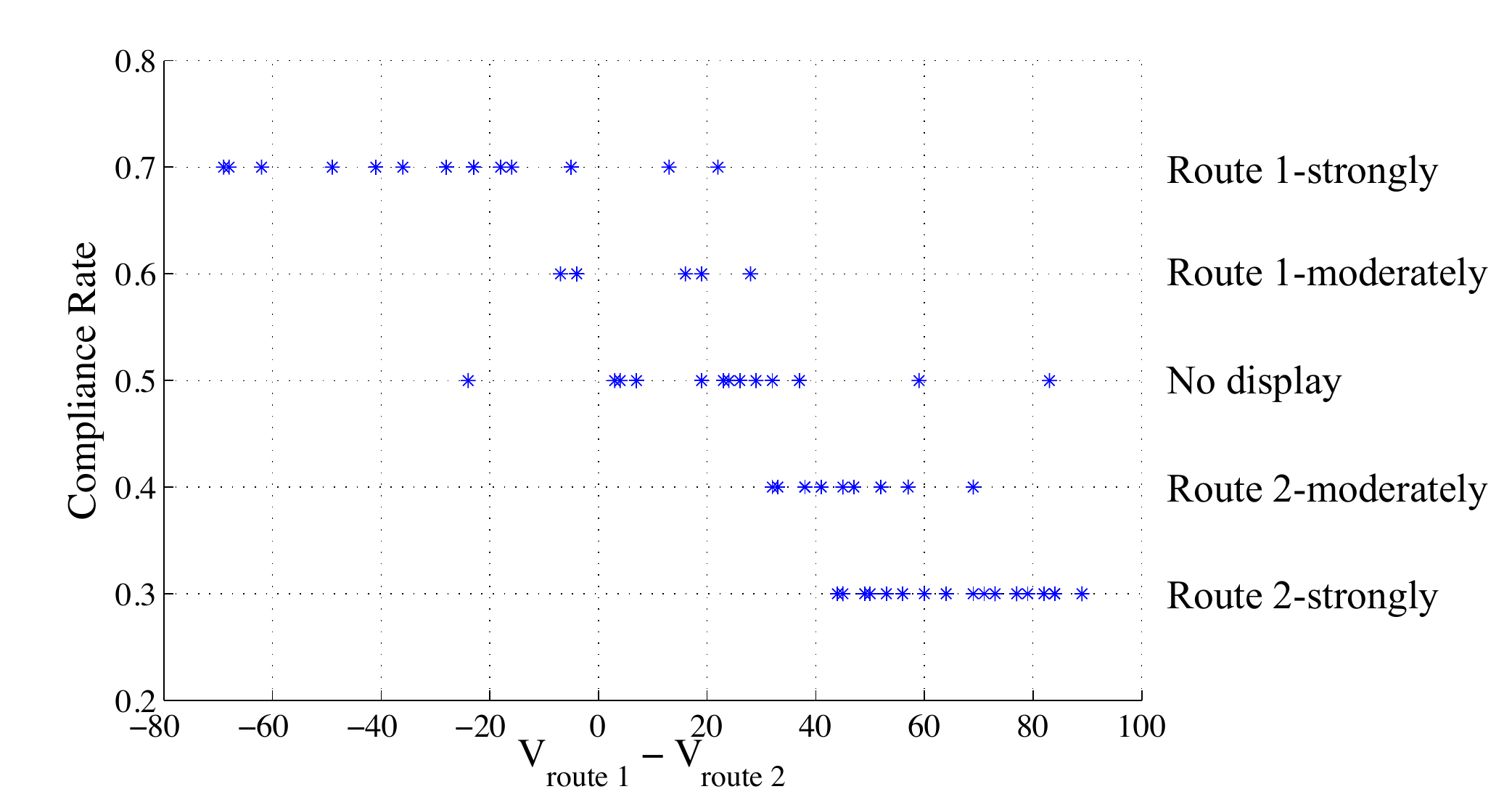}}
\caption{Comparison of the VMS message with the actual route condition. $V_{route 1}$ and $V_{route 2}$ are the traffic volume on Route 1 and Route 2, respectively.}
\label{Fig:VMSCompare}
\end{figure}

One of the key hypotheses of this study is the suspicion that completely factual and genuine display of route information may not be optimal on a system level, and a more sophisticated VMS strategy may sometimes presents routing suggestions that are not fully reflective of the actual traffic condition. This has been confirmed in Figure \ref{Fig:VMSCompare}, where we examine an optimal, LDR-based VMS strategy (with average travel time of 171.6 s where its genuine counterpart yields 178.3s). We match the displayed message with the actual traffic condition, represented as the difference of the traffic volume on the two alternative routes $V_{route 1}- V_{route 2}$. It is understood that when $V_{route 1}- V_{route 2}$ is very large, the VMS would recommend drivers to choose route $2$. However, although Figure \ref{Fig:VMSCompare} shows such a general trend overall, there are significant variations. For example, ``Route 1-strongly" is made for $V_{route 1}- V_{route 2}$ ranging from $-70$ to $20$, which seems counter-intuitive. Also, no display is made for some time instances when $V_{route 1}- V_{route 2}$ ranges between $-20$ and $80$. Overall, the VMS display seems to favor Route 1, and only recommends Route 2 when the traffic on Route 1 is significantly larger.

\section{Conclusion and future research}\label{secConclusion}
This paper explores the full potential of an on-line adaptive VMS display strategy for mitigating network-wide congestion. Unlike conventional display strategy, which tends to show accurate and factual information regarding the routing alternatives, we design the responsive VMS display in a more sophisticated way aiming at minimizing the network-wide delay while remaining robust and immune to network paradoxes that could arise in a complex traffic situation. 

In order to achieve this goal, we invoke the linear decision rule (LDR) approach for on-line traffic management. The LDR proposes a linear transformation from the current system state to the control parameters, which can be realized in real time. Then, the coefficients of such a transformation is optimized via the off-line module which boils down to a stochastic optimization problem. Implementation details of the LDR for the VMS control are provided in view of the real-world test network. The simulation-based test shows significant improvement of the LDR-based VMS display strategy over genuine route guidance. 

Future work will focus on: (1) the implementation and assessment of the proposed VMS framework on more realistic simulation platforms such as microsimulation supported by realistic data; (2) calibration of the modeling parameters, in particular the compliance rate which has a significant impact on the effectiveness of the VMS; (3) decision support for transportation authorities to devise effective route guidance.

\section*{Acknowledgement}
This work was supported by National Basic Research Program of China (973 Project) 2012CB725405, the National Science and Technology Support Program(2014BAG03B01), National Natural Science Foundation of China 61273238, Beijing Municipal Science and Technology Project (D15110900280000) and Tsinghua University Project (20131089307).

\bibliography{Reference}

\end{document}